\documentclass{article}

\usepackage{amsmath,amsthm,amsfonts,amssymb}

\theoremstyle{plain}
\newtheorem{theorem}{Theorem}[section]
\newtheorem{lemma}[theorem]{Lemma}

\newtheorem{proposition}[theorem]{Proposition}

\theoremstyle{definition}
\newtheorem{definition}[theorem]{Definition}

\theoremstyle{remark}

\newtheorem{remark}[theorem]{Remark}
\newtheorem{example}[theorem]{Example}

\newcommand{\tensor}{\ensuremath{\otimes}}

\usepackage{mathrsfs}
\usepackage{amsfonts}
\usepackage[all]{xy}
\usepackage{fancyhdr}
\usepackage[left=1in,right=1in,top=1in,bottom=1in]{geometry}

\def \C              {\mathcal C}

\def \I              {\mathcal I}

\def \CAA            {{\mathcal C}_{A/A}}

\def \M              {\mathcal M}
\def \NA             {{\mathcal N}_A}
\def \P              {\mathbb P}
\def \Po             {\mathbb{P}}
\def \U              {\mathcal U}

\begin{document}

\title{HKR Theorem for Smooth $S$-algebras}
\author{Randy McCarthy and Vahagn Minasian}
\date{}

\maketitle

\begin{abstract}
\noindent
We derive an \'etale descent formula for topological Hochschild homology 
and prove a HKR theorem for smooth $S$-algebras. 

\vspace{10pt}

\noindent
MCS: 55P42
\end{abstract}

\section{Introduction}

One of the main results for computing the Hochschild homology of smooth discrete 
algebras is the Hochschild-Kostant-Rosenberg (HKR) theorem (e.g. see Chapter 3 of 
~\cite{Loday}), which states that for a smooth algebra $k \to A$, the Hochschild
homology coincides with differential forms:
$$ HH_\ast (A) = \Omega^\ast _{A|k}.$$
In fact this result is often used not only to compute the Hochschild homology, but 
also the other way around: in order to generalize some results to non-smooth (or 
even non-commutative) algebras one replaces the differential forms by Hochschild 
homology. Other applications include a comparison theorem between cyclic and de Rham 
homology theories. 

One of our objectives is to develop a topological analogue of the HKR theorem in 
the framework provided in ~\cite{EKMM}, or more precisely, in the category of 
commutative of $S$-algebras. Recall that $S$-algebras are equivalent to the more 
traditional notion of $E_\infty$-ring spectra, and are a generalization to stable 
homotopy theory of the algebraic notion of a commutative ring.  
In this context, the topological Andr\'e-Quillen homology of a commutative 
$S$-algebra $A$ is the natural replacement 
of the module of differentials $\Omega^1 _{A|k}$, as it is evident from the definition 
of TAQ. The definitions of TAQ, as well as THH, in our context are recalled in Section 2,
and we refer to ~\cite{Maria} and Chapter IX of ~\cite{EKMM} for detailed discussion 
of these notions.  
 Noting that the orbits of 
the $n'th$ smash powers of the suspension 
of $TAQ$, $(\Sigma TAQ(A))^{\wedge_A n}/\Sigma_n$, are analogous to symmetric
powers in the graded context, and therefore correspond to 
taking exterior powers (and thus are the analogues of the higher order modules
of differentials), we state our main theorem.

\begin{theorem}
\label{th:main}
(HKR)
For a connective smooth $S$-algebra $A$, the natural (derivative) map 
$THH(A) \to \Sigma TAQ(A)$ has a section in the category of $A$-modules which 
induces an equivalence of $A$-algebras:

$$\P_A \Sigma TAQ(A) \stackrel{\simeq}{\to} THH(A),$$
where $\P$ is the symmetric algebra triple.  
\end{theorem}

The following is a description of the structure of the paper.

In Section 2 we recall the definitions of topological Hochschild homology and topological 
Andr\'e-Quillen homology in our framework. More precisely, the two main categories 
where our work takes place are the following. 
The first one is the 
category of $A$-modules, denoted by $\M_A$, where $A$ is a commutative $S$-algebra. 
There is a triple $\P_A : \M_A \to \M_A$ on this category given by 
$\Po M =\bigvee_{j \ge 0} M^j/ \Sigma_j$ (here $M^j$ denotes the $j$-fold smash 
power over $A$ and   $M^0 = A$), which leads us to the second category of interest - 
the category $\M_A[\P]$ of algebras in $\M_A$ over $\P$. Clearly, it is equivalent to 
the category of commutative $A$-algebras $\C_A$. For convenience, we denote the 
reduced version of $\P$ by $\P^1$. In other words, $\P^1$ is the obvious functor for 
which $\P = A \vee \P^1$. 

Note that both of these categories are closed model categories, and for a discussion 
on their homotopy categories we refer to Chapter VII of ~\cite{EKMM}. A good account 
for the general theory of closed model structures can be found in ~\cite{Dwyer}.

In Sections 3 and 4 we define \'etale, thh-\'etale, smooth and thh-smooth $S$-algebras, 
show that all these are generalizations of appropriate notions from discrete algebra, and 
prove their basic properties.   
 
Section 5 is devoted to establishing some conditions on a simplicial set $X_\ast$ 
and a map of commutative 
$R$-algebras $A \to B$ that imply the identity
\begin{equation}
\label{eq:1descent}
A \otimes X_\ast \wedge_A B \simeq B \otimes X_\ast.
\end{equation}
Observe that as a special case of this equation (more precisely, when we take the 
simplicial set $X_\ast$ to be the circle $S^1 _\ast$), we get an equation

\begin{equation}
\label{eq:2descent}
THH(A) \wedge_A B \simeq THH(B).
\end{equation}
Here we employed  the identity $THH(A) \simeq A \otimes S^1 _\ast$ derived by McClure,
Schw\"anzl and Vogt in ~\cite{Vogt}. Of course, in discrete algebra, the analogue of 
(~\ref{eq:2descent}) is referred to as {\it \'etale descent formula for HH} (see e.g.
~\cite{WG}). Following this, we will refer to both (~\ref{eq:1descent}) and 
(~\ref{eq:2descent}) as \'etale descent formulas.

To prove (~\ref{eq:1descent}), we produce a necessary condition for it to hold, and 
show that under some additional hypothesis, that condition is also sufficient. The 
notion of completeness is also discussed here, as it plays an important role in 
understanding  (~\ref{eq:1descent}).   
The equation (~\ref{eq:2descent}) is a key technical step in the proof of the HKR 
theorem for smooth $S$-algebras.

In Section 6, we prove the main (HKR) Theorem~\ref{th:main}, and conclude the section by  
showing that, as a consequence of the HKR theorem, 
the first fundamental sequence of the modules of differentials splits under a 
smoothness hypothesis. Here, following the terminology of discrete algebra, by 
first fundamental sequence of the modules of differentials we mean the homotopy 
cofibration sequence
$$TAQ^R(A) \wedge_A B \to TAQ^R(B) \to TAQ^A(B),$$
associated to the sequence $R \to A \to B$ of $S$-algebras (see ~\cite{Maria} for a
detailed discussion on this).   

Our definition for a map of commutative ring spectra $f:C\rightarrow D$
being \'etale when $TAQ(D|C)\simeq *$ is not new. We were first
introduced to this idea by F. Waldhausen in 1991.
Some other people whom we are aware of using this idea (either
formally or in private conversation) are:
M. Basterra,
T. Goodwillie,
T. Hunter,
J. Klein,
I. Kriz,
M. Mandell,
J.~McClure,
T. Pirashvili,
C. Rezk,
B. Richter,
A. Robinson,
J. Rognes,
J. Smith,
and
S.~Whitehouse. 
The idea of thh-\'etale that we use seems fairly common
to the extent that 
most of these people have considered this also.  
In particular, recent work by  J. Rognes independently
establishes several of the structural properties of
thh-\'etale maps which we use. 

We have been greatly aided by many mathematicians while
working out our ideas for this paper. In particular, we
would like to thank Maria Basterra for teaching us about
commutative $S$--algebras and how to work with them.
We thank Mike Mandell for his support, insights and important examples.
This work arose from a series of talks with Charles Rezk (who
also caught a serious mistake in an earlier draft)  while he
taught us about the DeRham cohomology of commutative ring 
spectra. We came upon the main conjecture while talking
with Birgit Richter and were certainly motivated by ideas of Nick
Kuhn about splitting Goodwillie Taylor towers.

\section{Preliminaries: THH and TAQ of commutative $S$-algebras}

In this section we give a brief introduction into THH and TAQ of commutative 
$S$-algebras.
Chapter IX of ~\cite{EKMM} and ~\cite{Maria} provide a good in depth discussion 
of these notions in our framework. 

Let $R$ be a cofibrant commutative $S$-algebra, $A$ - a cofibrant $R$-algebra or 
a cofibrant commutative $R$-algebra, and $M$ an $(A,A)$-bimodule. Write $A^p$ for 
the $p$-fold $\wedge_R$-power, and let 
$$\phi:A \wedge_R A \to A \hspace{.5in} and \hspace{.5in} \eta:R \to A$$
be the product and unit of $A$-respectively.

Let 
$$\xi_l: A \wedge_R M \to M \hspace{.5in} and \hspace{.5in} \xi_r:M \wedge_R A \to M$$
be the left and right actions of $A$ on $M$. Denote the canonical cyclic permutation
isomorphism by $\tau$:
$$\tau:M \wedge_R A^p \wedge_R A \to A \wedge_R M \wedge_R A^p.$$ 

\begin{definition}
Let $THH^R(A;M)_\ast$ be the simplicial $R$-module whose $R$-module of $p$-simplices
is $M \wedge_R A^p$, and whose face and degeneracy operators are 
$$d_i = \left\{ \begin{array}{clcr}
           \hspace{-1in} \xi_r \wedge (id)^{p-1}               & if & i=0 \\
           id \wedge(id)^{i-1} \wedge \phi \wedge (id)^{p-i-1} & if & 1 \leq i <p\\
           \hspace{-.7in}(\xi_l \wedge (id)^{p-1}) \circ \tau  & if & i=p
          \end{array}
  \right .$$
$$ s_i= id \wedge (id)^i \wedge \eta \wedge (id)^{p-i}.$$  
Define 
$$THH^R(A;M) = |THH^R(A;M)_\ast|.$$
When $M=A$, we delete it from the notation, writing $THH^R(A)$.

\end{definition}

Clearly this definition(~\cite{EKMM})mimics the definition of the standard complex
for the computation of Hochschild homology, as given in ~\cite{Cart}. Of course, the 
passage from a simplicial spectrum to its geometric realization is the topological 
analogue of passage from a simplicial $k$-module to a chain complex. 

Observe that the maps
$$\xi_p=id \wedge \eta^p:M \simeq M \wedge_R R^p \to M \wedge_R A^p$$
induce a natural map of $R$-modules
$$\xi=|\xi_\ast|:M \to THH^R(A;M).$$
If $A$ is a commutative $R$-algebra, then clearly $THH^R(A)_\ast$ is a simplicial
commutative $R$-algebra and $THH^R(A;M)_\ast$ is a simplicial $THH^R(A)$-module. 
Hence, $THH^R(A)$ is a commutative $A$-algebra with the unit map given by the above map
$\xi:A \to THH^R(A)$. 

Observe that if $M$ is an $(A,A)$-bimodule and $\underline{M}$ is the corresponding 
constant simplicial $(A,A)$-bimodule, then
$$M \wedge_{A^e}\beta^R(A) \cong|\underline{M}\wedge_{A^e}\beta^R_\ast (A)|,$$
where $A^e=A \wedge A^{op}$. We have canonical isomorphisms
$$M \wedge_R A^p \cong M \wedge_{A^e}(A^e \wedge_R A^p) \cong M 
    \wedge_{A^e}(A \wedge_R A^p \wedge_R A)$$
given by permuting $A^{op}=A$ past $A^p$. As these isomorphism commute with the face
and degeneracy operations, we get
\begin{equation}
THH^R(A;M) \cong M \wedge_{A^e} \beta^R(A). \label{eq:barthh}
\end{equation}

Now we turn our attention to the Topological Andr\'e-Quillen Homology. The definition,
presented by Maria Basterra in ~\cite{Maria}, employs the following two functors. 

{\bf The augmentation ideal functor.}
Let $A$ be a commutative $S$-algebra, and $I: \CAA \to \NA$ the functor from the category 
of commutative $A$-algebras over $A$ to the category of $A$-NUCA's which assigns to 
each algebra ($B$, $\eta:A \to B$, $\epsilon:B \to A$) its ``augmentation ideal'': $I(B)$
defined by the pullback diagram in $\M_A$,
$$
\xymatrix{
I(B) \ar[r] \ar[d]
&
B \ar[d]^\epsilon\\
\ast \ar[r]
&
A.
}$$  
Note that by the universal property of pullbacks $I(B)$ comes with a commutative 
associative (not necessarily unital) multiplication. Moreover, this functor has a 
left adjoint $K:\NA \to \CAA$ which maps a non-unital algebra $N$ to $N \vee A$
(Proposition 3.1 of ~\cite{Maria}). In addition this adjunction produces an
equivalence of homotopy categories given by the total derived functors 
${\bf L}K$ and ${\bf R}I$ (Proposition 3.2 of ~\cite{Maria}).

{\bf The indecomposables functor.}
Let $Q:\NA \to \M_A$ denote the ``indecomposables'' functor that assigns to each $N$ in 
$\NA$ the $A$-module $Q(N)$ given by the pushout diagram in $\M_A$
$$
\xymatrix{
N \wedge_A N \ar[r] \ar[d]
&
\ast \ar[d]\\
N \ar[r]
&
Q(N).
}$$   
This functor has a right adjoint $Z:\M_A \to \NA$ given by considering $A$-modules as 
non-unital algebras with zero multiplication. Since $Z$ is the identity on morphisms
and the closed model structure on $\NA$ is created in $\M_A$, $Z$ preserves fibrations 
and acyclic fibrations, so by Chapter 9 of ~\cite{Dwyer}, the total derived functors
${\bf R}Z$ and ${\bf L}Q$ exist and are adjoint.  

\begin{definition}

Let $B \to A$ be a map of commutative $S$-algebras. Define 
$$ TAQ(B/A)=\Omega_{B/A} \stackrel{def}{=} {\bf L}Q {\bf R}I (B \wedge^{\bf L} _A B),$$
where $B \wedge^{\bf L} _A B$ denotes the total derived functor of 
$- \wedge^{\bf L} _A B$ evaluated at $B$.

\end{definition}

Of course, as it is observed in ~\cite{Maria}, $\Omega_{B/A}$ is simply a derived 
analogue of the $B$-module of K\"ahler differentials from classical algebra. 

{\bf Notation.}
Fix a cofibrant commutative $S$-algebra $A$. Then for an $A$-algebra $B$ and an 
$A$-module $M$, we denote by $THH(B,M|A)$ and $TAQ(B,M|A)$ the topological Hochschild
and Andr\'e-Quillen homologies of $B$ over $A$ with coefficients in $M$. 
If $M=B$, we omit it from the notation. In addition, $\widetilde{THH}(B|A)$ stands
for the reduced topological Hochschild homology, defined to be the homotopy cofiber of 
the natural map $B \to THH(B|A)$.

\section{(thh-)\'etale $S$-algebras}

Recall that in discrete algebra smooth maps can be roughly defined to be the maps 
which can be decomposed into a polynomial extension followed by an \'etale extension. 

\begin{definition}
We say that a discrete $k$-algebra $A$ is smooth if for any prime ideal of $A$ there 
is an element $f$ not in that prime such that there exists a factorization
$$k \to k[x_1, \cdots ,x_n] \stackrel{\phi}{\to} A_f$$
with $\phi$ \'etale, i.e. flat and unramified. 
\end{definition}

Under some finiteness and flatness conditions this notion of smooth maps coincides 
with most other standard ones (see Appendix of ~\cite{Loday}). It is with this approach
to smoothness in mind that we define our smooth maps of $S$-algebras. Hence the need 
to discuss the notion of \'etale algebras first. Recall that for discrete algebras, both 
smooth and \'etale maps are defined to be finite in some appropriate sense. We do not 
impose a finiteness condition on $S$-algebras as it is not needed for our main results.
Consequently, a more appropriate terminology to use would be `formally' \'etale and 
smooth, which we don't for the sake of economy.    

We begin with a pair of definitions.
Let $R$ be a commutative cell $S$-algebra and $A$, $C$ and $D$  commutative  
$R$-algebras.

\begin{definition}
The map of algebras $C \to D$ is \'etale (thh-\'etale) if $TAQ(D|C)$ is contractible
($D \stackrel{\simeq}{\to} THH(D|C)$).
\end{definition}

We also define (thh-)\'etale coverings to be  faithfully flat families of 
(thh-)\'etale extensions: 

\begin{definition}
We say that $\{ A \to A_\alpha \}_{\alpha \in \I}$ is a (thh-)\'etale covering of $A$ 
if

1. each map  $ A \to A_\alpha $ is (thh-)\'etale, and 

2. for each pair of $A$-modules $M \to N$ such that 
$M \wedge A_\alpha \to N \wedge A_\alpha$ is a weak equivalence for all $\alpha$, 
the map $M \to N$ is itself a weak equivalence. 
\end{definition}

This definition gives rise to a few natural questions. Are there `enough' (thh-)\'etale
coverings? What is the relationship between \'etale and thh-\'etale?

\begin{remark}
We claim that for each commutative $R$-algebra $A$, at least one (non-trivial)
\'etale covering and one (non-trivial) thh-\'etale covering exists. 
To see this, first recall some facts about localizing $S$-algebras.  

Suppose $T$ is a 
multiplicatively closed subset of $\pi _\ast (A)$. Then by Section 1 of Chapter V of 
~\cite{EKMM}, for each $A$-module $M$ one can define a localization $M[T^{-1}]$ of 
$M$ at $T$ using a telescope construction with a key property

\begin{equation}
\pi _\ast   (M[T^{-1}]) \cong \pi _\ast (M)[T^{-1}].
\end{equation}
Moreover, the localization of $M$ is the smash product of $M$ with the 
localization of $A$. 
In addition, by Theorem VIII.2.1 of ~\cite{EKMM} one can construct the 
localization in such a 
way that $A[T^{-1}]$ is a cell $R$-algebra and the localization map 
$A \to A[T^{-1}]$ is an inclusion of a subcomplex. 
Moreover, since $A[T^{-1}]$ smashed over $A$ with itself is equivalent to localizing 
$A[T^{-1}]$ at $T$, we conclude that $A[T^{-1}] \wedge_A A[T^{-1}] \cong A[T^{-1}]$, 
and hence the map $A \to A[T^{-1}]$ is (thh)-\'etale. 
Now for each prime of $\pi_\ast (A)$, pick an element $f$ outside of it, and let $T$
be the multiplicative system generated by that element. Let $M \to N$ be a map 
of $A$-modules such that $M_f \to N_f$ is an equivalence for all $f$. In other words,
the induced map $\pi_\ast(M) \to \pi_\ast (N)$ is such that the localizations of this map 
are isomorphisms. Hence the map itself is an isomorphism (e.g. see Chapter 2 of 
~\cite{Eisen}), proving that $\{A \to A_f\}$ is a covering. Of course, there are other 
collections of multiplicative systems in $\pi_\ast (A)$ that we can use to produce a 
covering (e.g. all the maximal ideals of $\pi_\ast (A)$); the key property is that if 
a map of modules localized at these systems is an isomorphism then the map itself is an 
isomorphism.      

\end{remark}

Recall that the Goodwillie derivative of $THH$ is the suspension of $TAQ$ and thus 
thh-\'etale implies \'etale. 
This is discussed in detail for example in ~\cite{Minas}. While 
the converse is false in general, it does hold for certain classes of spectra; for 
example, the two notions are equivalent for connective spectra (see~\cite{Minas}).   
The following example (communicated by M. Mandell, ~\cite{Mandell}) illustrates that 
\'etale does not always imply thh-\'etale.  

\begin{example}
We work over the field ${\bf{F}}_p$. Fix $n>1$ and let 
$C^\ast (K({\bf Z}/p{\bf Z}, n))$ be 
the cochain complex of $K({\bf Z}/p{\bf Z}, n)$ viewed as an $E_\infty$-algebra. 
To ease the notation we denote this $E_\infty$-algebra by $R$. 
$R$ has a non-zero homotopy group in degree $-n$, while its $-n+1$'st homotopy group
is trivial. Recall that $THH(R|{\bf{F}}_p)$ is equivalent to $Tor^{R \otimes R}(R,R)$, 
hence we have an Eilenberg-Moore type spectral sequence (see Theorem IV.6.2 or 
Theorem IX.1.9 of ~\cite{EKMM}):
$$Tor_{p,q} ^{\pi_\ast(R \otimes R)} ( \pi_\ast (R),\pi_\ast (R))
\Rightarrow  Tor_{p+q} ^{R \otimes R}(R,R)= THH_{p+q}(R|{\bf{F}}_p).$$
Consequently, the $-n+1$'st homotopy group of $THH(R|{\bf{F}}_p)$ is non-trivial.
Hence $R$ and $THH(R|{\bf{F}}_p)$ are not equivalent, and thus, $R$ is not thh-\'etale. 

To see that $R$ is \'etale we need to give another description for $R$ that requires 
the use of generalized Steenrod operations for $E_\infty$-algebras 
(see ~\cite{May1} for a reference on Steenrod operations in our context).
In fact, we will only need the operation $P^0$. Recall that it preserves degree 
and performs the $p$'th power operation on elements in degree 0. 
By Section 6 of~\cite{Mike}, $R$ can be described as the 
$E_\infty$-algebra free on two generators $x$ (in degree $-n$) and $y$
with $dx=0$ and $dy=x- P^0 x$. Then noting that $P^0 x$ is of the form 
$e \otimes x^{\otimes p}$, where $e$ is in $E(p)$ ($E$ being the $E_\infty$ operad), 
we observe that the $R$-module representing $TAQ(R)$ is modeled by the free $R$-module
on two generators $\bar{x}$ and $\bar{y}$ with 

$d \bar{x} =0$ and

$d  \bar{y} = \bar{x} -
e \otimes[\bar{x} \otimes x \otimes \cdots \otimes x +
          \cdots x \otimes \cdots \otimes x \otimes \bar{x}]
= \bar{x} - e (1 +a + \cdots + a^{p-1}) \otimes
            [\bar{x} \otimes x \otimes \cdots \otimes x]$,

\noindent
where $a$ is a generator of of the cyclic group of $p$ elements. Observe that we 
have an $R$-module contraction $s$ given by

$s(\bar{y}) = 0$ and

$s(\bar{x}) = \bar{y}+f \otimes [\bar{x} \otimes x \otimes \cdots \otimes x]$,

\noindent
where $f$ is such that $df = e (1 +a + \cdots + a^{p-1})$. Thus $TAQ(R)$ is contractible. 

\end{example}

In the following lemma we prove a few easy properties of \'etale maps that will be needed
later. 

\begin{lemma}
\label{lemma:etale}
1. (Transitivity) If $A$ is \'etale over $R$ and $B$ is \'etale over $A$, 
then $B$ is \'etale over $R$. 

2. (Base Change) If $B$ and $C$ are cofibrant $A$ algebras and $B$ is \'etale over 
$A$, then $C \wedge_A B$ is \'etale over $C$. Also, if $B \to C$ is a 
\'etale map of $A$-algebras and $D$ is a cofibrant $A$-algebra then 
$B \wedge_A D \to C \wedge_A D $ is also \'etale. 

3. (Polynomial Extensions) If $B \to C$ is a \'etale map of $A$-algebras, then 
for all cell $A$-modules $X$, the induced map 
$\P _B (X \wedge_A B) \to  \P _C (X \wedge_A C)$ is also \'etale.
\end{lemma}

While the lemma and the following proof are stated for \'etale extensions, a similar
result holds for thh-\'etale algebras as well. 
The remark after the lemma describes how to 
adjust the proof for the thh-\'etale case. 

\begin{proof}

1. The transitivity is immediate from the cofibration sequence induced by 
$R \to A \to B$:
$$TAQ(A|R) \wedge_A B \to TAQ(B|R) \to TAQ(B|A).$$

2. By Proposition 4.6 of ~\cite{Maria}, 
$TAQ(C \wedge_A B |C) \simeq TAQ(B|A) \wedge_A C$. Since $TAQ(B|A) \simeq \ast$, 
$TAQ(C \wedge_A B |C)$ is also contractible. Now let  $B \to C$ be an \'etale 
map, then for any $A$-algebra $D$, 
$$TAQ(C \wedge_A D|B \wedge_A D) \simeq
TAQ(C \wedge_B B \wedge_A D|B \wedge_A D) \simeq TAQ(C|B) \wedge_B  B \wedge_A D.$$
Here the second map is an equivalence by Proposition 4.6 of ~\cite{Maria} once again. 
Recalling that the map $C \to B$ is \'etale, we conclude that 
$TAQ(C \wedge_A D|B \wedge_A D) \simeq \ast$. 

3. It is immediate from part 2, once we observe that 
$\P _B (X \wedge_A B) \cong \P _A (X) \wedge_A B$. 

\end{proof}

\begin{remark}
\label{remark:thhtaq}
Note that the proof of Lemma~\ref{lemma:etale} (\'etale case) hinges on two key facts about 
$TAQ$:

1.For cofibrant $A$ algebras $B$ and $C$, $TAQ(C \wedge_A B |C) \simeq TAQ(B|A) \wedge_A C$.

2. If the map of $A$-algebras $C \to D$ is \'etale then 
$TAQ(C|A) \wedge_C D \simeq TAQ(D|A)$. 

Thus, if analogous results hold for $THH$, then the arguments of the above proof can be 
repeated to prove the lemma in the thh-\'etale case. In fact, this reasoning also extends to 
future results (e.g. Lemma~\ref{lemma:smooth}), in which the \'etale assumption may be 
replaced by the thh-\'etale one.   

To see the analogue of the first fact about $THH$, simply recall the definition of $THH$
that mimics the standard complex for the computation of algebraic Hochschild homology 
(see ~\cite{EKMM}). Then  $THH(C \wedge_A B |C)$ and $THH(B|A) \wedge_A C$ both have 
$B \wedge_A \cdots \wedge_A B \wedge_A C$ as simplices and the map between them is the 
identity map on simplicial level. Thus the two objects are equivalent. 

The analogue of the second fact (with some extra conditions) is listed as 
Lemmas~\ref{etdescomp} and ~\ref{etdesconn} and will be proved 
later. 
\end{remark}

We have the following result about \'etale maps.

\begin{proposition}
\label{prop:emetale}
1. If $A \to B$ is an \'etale map of discrete algebras, then the induced 
map of $S$-algebras $HA \to HB$ is also \'etale. 

2. If for a commutative ring $k$, $h:Hk \to B$ is an \'etale map of $S$-algebras, 
then the map $Hk \to H \pi_0(B)$ which realizes the map induced by $h$ on 
$\pi_0$ is also \'etale.
\end{proposition}

\begin{proof}
1. Let $A \to B$ is an \'etale map of discrete algebras. We need to show that 
$TAQ(HB|HA)$ is contractible. Since $HA$ and $HB$ are connective this is equivalent
to showing that the natural map $\phi:HB \to THH(HB|HA)$ is a weak equivalence. 

Since $A \to B$ is \'etale, it is in particular flat, hence by 
Theorem IX.1.7 of ~\cite{EKMM}, $\pi_\ast (THH(HB|HA)) \cong HH_\ast (B|A)$. 
However for \'etale maps we have that $HH_0(B|A) \cong B$ and 
$HH_\ast (B|A) \cong \Omega^\ast _{B|A} = 0$ for $\ast>0$. 
Thus, $\phi$ induces an isomorphism
on $\pi_\ast$ for $\ast>0$ as it is simply the unique map between trivial groups. 
Combining this with the fact that $\phi$ on $\pi_0$ is the identity map on $B$, we 
conclude that $\phi$ is a weak equivalence. 

2. Let $Hk \to B$ be an \'etale map of $S$-algebras. $B$ is a generalized 
Eilenberg-MacLane spectrum since it is a module over the Eilenberg-MacLane spectrum
$Hk$. Hence there is a map $f:H \pi_0 (B) \to B$ that realizes the identity 
map on $\pi_0$. Also, for the same reason, we have a map of 
$Hk$-algebras $g:B \to H \pi_0 (B)$ that induces the identity map on $\pi_0 (B)$. 
The sequence $H \pi_0 (B) \to B \to H \pi_0 (B)$ produces a pair of maps:
\begin{equation}
\label{eq:eilmac}
TAQ(H \pi_0 (B)|Hk) \to TAQ(B|Hk) \to TAQ(H \pi_0 (B)|Hk).
\end{equation}
Since $g \circ f$ is the identity, the composite map (~\ref{eq:eilmac}) is 
also an equivalence. However, $TAQ(B|Hk) \simeq \ast$, since $Hk \to B$ is \'etale. 
Hence $TAQ(H \pi_0 (B)|Hk) \simeq \ast$, proving that $Hk \to H \pi_0 (B)$ is \'etale.   


\end{proof}

We already mentioned that localizations provide a large class of examples of 
(thh-)\'etale maps. As in discrete algebra, another principal source of examples is
given by Galois extensions. The following definition is due to John Rognes 
(~\cite{Rognes}).

\begin{definition}
Let $B$ be a cofibrant $A$-algebra, and
$G$ be a grouplike topological monoid acting on $B$ through $A$-algebra maps, 
such that $G \simeq \pi_0(G)$ is finite. Then $A \to B$ is a $G$-Galois extension 
if 

(1) $A \simeq B^{hG}=F(EG_+,B)^G$, and

(2) $B \wedge_A B \simeq F(G_+,B)$,

\noindent
where $F$ is the internal function spectrum (see Section I.7 of ~\cite{EKMM}).
\end{definition}

\begin{proposition}
(Rognes)
A $G$-Galois extension $A \to B$ is thh-\'etale (and hence also \'etale).
\end{proposition}

\begin{proof}
$B \wedge_A B \simeq F(G_+,B)$ is a product of copies of $B$ so $B$ is a retract 
of $B \wedge_A B$. Hence the composite 
$B \to THH(B|A)=THH(B,B|A) \to THH(B, B \wedge_A B|A) \simeq B$ is an equivalence
and the last map splits (via the retract map). Moreover, since  $B \wedge_A B$ is 
a product of copies of $B$ that map is also a monomorphism in the derived category.
Hence, $B \to THH(B|A)$ is an equivalence.
\end{proof}

For examples of Galois extension we again refer to ~\cite{Rognes}.

\section{Smooth $S$-algebras}

\begin{definition}
The map of algebras $f:R \to A$ is (thh-)smooth if there is a (thh-)\'etale covering 
$\{ A \to A_\alpha \}_{\alpha \in \I}$ of $A$ such that for each $\alpha$ there is 
a factorization
$$ R \longrightarrow \P_R X \stackrel{\phi}{\longrightarrow} A_\alpha, $$
where $X$ is a cell $R$-module and $\P_R X$ is the free commutative $R$-algebra 
generated by $X$, with $\phi$ (thh-)\'etale. 
\end{definition}

As always, we would like the smooth $S$-algebras to generalize the corresponding
notion from discrete algebra. Let $k \to A$ be a smooth map of discrete algebras, 
in other words, for each prime ideal of $A$, there is an element $f$ away from it 
such that there is a factorization 
$k \to k[x_1, \cdots , x_n] \stackrel{\phi}{\to} A_f$ with $\phi$ \'etale.
We claim that $Hk \to HA$ is a smooth map of $S$-algebras. Indeed, we have a pair 
of maps $Hk \to Hk[x_1, \cdots , x_n] \stackrel{H\phi}{\to} HA_f$, where $H\phi$ is 
\'etale by Proposition~\ref{prop:emetale}. By the same proposition, we also get 
that $HA \to HA_f$ is \'etale. Moreover, the maps $HA \to HA_f$ form a covering, 
as smashing with $HA_f$ over $HA$ is equivalent to localizing at $f$. Thus, 
observing that $Hk[x_1, \cdots , x_n] \cong \P_{Hk}(\bigvee_n Hk)$, we conclude 
that $Hk \to HA$ is smooth. 

In the following lemma we list some of the basic properties of (thh-)smooth 
$S$-algebras.
Before doing so, we recall that the localization at a cell $R$-module $E$ is 
called smashing if for all cell $R$-modules $M$, the localization of $M$ at $E$ 
is given by $R_E \wedge_R M$, where $R_E$ is the localization of $R$ at $E$. 

\begin{lemma}
\label{lemma:smooth}

1. (Localization) If $A$ is (thh-)smooth over $R$ and the localization at $E$ is 
smashing, then the composite map $R \to A_E$ is also (thh-)smooth. 

2. (Transitivity) If $A$ is (thh-)smooth over $R$ and $B$ is (thh-)smooth over $A$,
then $B$ is (thh-)smooth over $R$.

3. (Base Change) If $A$ is (thh-)smooth over $R$, and $R \to B$ is a map of 
commutative $S$-algebras, then $B \to A \wedge_R B$ is also (thh-)smooth. 

\end{lemma}

\begin{proof}
Again, we present a proof of the smooth case. As noted in Remark~\ref{remark:thhtaq},
the proof of thh-smooth case is identical to this one.  

1. Since the localization at $E$ is smashing, $A_E \wedge_A A_E$ is the localization 
of $A_E$ at $E$. However, $A_E$ is already $E$-local. Hence the multiplication map
$A_E \wedge_A A_E \to A_E$ is an equivalence, implying that $TAQ(A_E|A) \simeq \ast$.
In other words, $A \to A_E$ is \'etale. Thus, it is smooth, since for the \'etale 
covering required by the definition of smoothness we can simply take the identity 
map of $A_E$. So the localization property will follow once we prove the transitivity
of smooth algebras. 

2. Let $A \to A_\alpha$ and $B \to B_\beta$
be \'etale coverings of $A$ and $B$ respectively such that there are
factorizations $ R \to \P_R(X_\alpha) \stackrel{\phi_\alpha}{\to} A_\alpha$ and 
$ A \to \P_A(Y_\beta) \stackrel{\psi_\beta}{\to} B_\beta$ with
$\phi_\alpha$ and $\psi_\beta$ \'etale. Consider the maps 
\begin{equation}
\label{eq:trans}
B \rightarrow B_\beta \wedge_A A_\alpha \wedge_{\P_R(X_\alpha)} A_\alpha .
\end{equation}
By parts 2 of Lemma~\ref{lemma:etale}, we have that the maps
$B \to B_\beta \wedge_A A_\alpha$ and
$A_\alpha \to A_\alpha \wedge_{\P_R(X_\alpha)} A_\alpha$ are \'etale. Hence,
the map $B_\beta \wedge_A A_\alpha \to 
B_\beta \wedge_A A_\alpha \wedge_{\P_R(X_\alpha)} A_\alpha$ is also \'etale.
Thus, by transitivity of \'etale extensions (part 1 of Lemma~\ref{lemma:etale}), 
we get that the above maps ~\ref{eq:trans} are \'etale. Next we show that this 
collection of \'etale maps forms a covering. To see this first observe that since
$A \to A_\alpha$ is a covering of $A$, so is
$A \to A_\alpha \wedge_{\P_R(X_\alpha)} A_\alpha$, as the multiplication map 
$A_\alpha \wedge_{\P_R(X_\alpha)} A_\alpha \to A_\alpha$ splits. 
Now let $M \to N$ be $B$-modules such that 
$M \wedge_B B_\beta \wedge_A A_\alpha \wedge_{\P_R(X_\alpha)} A_\alpha 
\stackrel{\simeq}{\to}
N \wedge_B B_\beta \wedge_A A_\alpha \wedge_{\P_R(X_\alpha)} A_\alpha$. 
Since $A \to A_\alpha \wedge_{\P_R(X_\alpha)} A_\alpha$ is a covering,
we conclude that for each $\beta$, 
$M \wedge_B B_\beta \stackrel{\simeq}{\to} N \wedge_B B_\beta$, and hence 
$M \simeq N$. 

Thus, it remains to show that 
$B_\beta \wedge_A A_\alpha \wedge_{\P_R(X_\alpha)} A_\alpha$ is \'etale over a 
polynomial extension of $R$. By part 2 of Lemma~\ref{lemma:etale} we have that 
$\P_A(Y_\beta)\wedge_A A_\alpha \wedge_{\P_R(X_\alpha)} A_\alpha \to
B_\beta \wedge_A A_\alpha \wedge_{\P_R(X_\alpha)} A_\alpha$ is \'etale. 
Now we simply observe that 
$\P_A(Y_\beta)\wedge_A A_\alpha \wedge_{\P_R(X_\alpha)} A_\alpha \cong
\P_{A_\alpha} (Y_\beta \wedge_A A_\alpha ) \wedge_{\P_R(X_\alpha)} A_\alpha \cong
\P_{\P_R(X_\alpha)}(Y_\beta \wedge_A A_\alpha \wedge_{\P_R(X_\alpha)} A_\alpha)$
and the last object being a polynomial extension of a polynomial over $R$ is itself
a polynomial over $R$.

3. Let $A \to A_\alpha$ be an \'etale covering of $A$ such that there are
factorizations $ R \to \P_R(X_\alpha) \stackrel{\phi_\alpha}{\to} A_\alpha$ with
$\phi_\alpha$ \'etale. For any $R$-algebra $B$,  by part 2 of Lemma ~\ref{lemma:etale}, 
the maps $B \wedge_R A \to B \wedge_R A_\alpha$ are \'etale. Moreover, since 
$\P_R(X_\alpha) \stackrel{\phi_\alpha}{\to} A_\alpha$ are \'etale, so are 
$\P_R(X_\alpha) \wedge_R B \to A_\alpha \wedge_R B$.
Note that $\P_R(X_\alpha) \wedge_R B \cong \P_B (X_\alpha \wedge_R B)$. Thus,
we have factorizations 
$B \to \P_B (X_\alpha \wedge_R B) \stackrel{\psi_\alpha}{\to} A_\alpha \wedge_R B$
with $\psi$ \'etale. 

To complete the proof it remains to show that the collection of \'etale maps
$B \wedge_R A \to B \wedge_R A_\alpha$ forms an \'etale covering. 
Let $M \to N$ be a pair of $B \wedge_R A$-modules such that 
$M \wedge_{B \wedge_R A} B \wedge_R A_\alpha  \stackrel{\simeq}{\to} 
N \wedge_{B \wedge_R A} B \wedge_R A_\alpha$. Observe that 
$$M \wedge_{B \wedge_R A} B \wedge_R A_\alpha \cong
M \wedge_{B \wedge_R A} B \wedge_R A \wedge_A A_\alpha \cong
M \wedge_A A_\alpha.$$ 
Thus we get that 
$ M \wedge_A A_\alpha \stackrel{\simeq}{\to} N \wedge_A A_\alpha$, and since 
$A \to A_\alpha$ is an \'etale covering, we conclude that $M \simeq N$.  

\end{proof}

\section{\'Etale Descent}

Our main goal is to prove the topological analogue of the HKR theorem. 
As will be observed later, it is of critical importance for HKR that we be able 
to identify conditions on the map of $R$-algebras $A \to B$, that will imply the 
identity $THH(A|R) \wedge_A B \simeq THH(B|R)$. In fact recalling that by ~\cite{Vogt}
$THH(A|R) \cong A \otimes_R S^1$, we can rewrite the above identity as 
$(A \otimes_R S^1) \wedge_A B \simeq B \otimes_R S^1$, which prompts us to investigate
conditions on a simplicial set $X$ and a map of $R$-algebras $A \to B$ that imply 
the more general identity

\begin{equation}
\label{eq:general}
(A \otimes_R X) \wedge_A B \simeq B \otimes_R X.
\end{equation}  
Almost immediately we can get a necessary condition for (~\ref{eq:general}) to hold. 
First we need a change of base formula for tensor products
\begin{equation}
\label{eq:base}
 A \wedge_{A \otimes_R X} (B \otimes_R X) \simeq B \otimes_A X.
\end{equation}
We are grateful to M. Mandell for suggesting a proof of this formula 
by describing the $A$-algebra maps into a fixed $A$-algebra $C$.

First consider $\C _A (A \wedge_{A \otimes_R X} (B \otimes_R X), C)$. By universal 
property of pushouts, this is isomorphic to the subset of maps $f$ in 
$\C _R (B \otimes_R X, C)$, such that the restriction of $f$ to $A \otimes_R X$ factors
through $A \otimes_R X \to A$. By adjunction of the tensor product, 
$\C _R (B \otimes_R X, C) \cong \U (X, \C _R (B,C))$. Thus, 
$\C _A (A \wedge_{A \otimes_R X} (B \otimes_R X), C)$ is isomorphic to the subset of maps 
$\phi$ in $\U (X, \C _R (B,C))$ such that for all $x \in X$, $\phi (x):B \to C$ restricted 
to $A$ is the same map, in other words, the maps $A \to B \stackrel{\phi (x)}{\to} C$ and
$A \to B \stackrel{\phi (y)}{\to} C$ are the same for all $x,y \in X$. Observe that the 
collection of such maps is precisely $\U (X, \C _A (B,C)) \cong \C _A (B \otimes_A X, C)$, 
and hence the proof of the formula (~\ref{eq:base}) is complete by Yoneda's lemma. 

Now consider the following commutative diagram of $A$-algebras

$$
\xymatrix{
A \wedge_{A \otimes_R X} (A \otimes_R X) \wedge_A B \ar[r] \ar[d]
&
A \wedge_{A \otimes_R X} (B \otimes_R X) \ar[d]\\
B \ar[r]
&
B \otimes_A X .
}
$$
The left vertical arrow is clearly an isomorphism, and by the base change formula
(~\ref{eq:base}), so is the right vertical arrow. 
Hence if we assume that the identity (~\ref{eq:general}) holds, then the top 
horizontal map is an equivalence, implying that the bottom map 
$B \to B \otimes_A X$ is also an equivalence. 

Thus, $B \stackrel{\simeq}{\to} B \otimes_A X$ is a necessary condition for 
(~\ref{eq:general}) to hold. Of course, in general this condition alone is not enough 
to ensure (~\ref{eq:general}), as can easily be seen on example of $X=S^0$. 
$B \otimes_A S^0 \cong B \wedge_A B$, and hence the condition 
$B \stackrel{\simeq}{\to} B \otimes_A X$ becomes 
$B \stackrel{\simeq}{\to} B\wedge_A B$. This, however, as can be seen in the following 
example, does not imply 
$A \wedge_R B \simeq B \wedge_R B$, which is the restatement of (~\ref{eq:general})
for $X=S^0$.

\begin{example}
Consider a pair of $S$-algebra maps 
$H{\bf Z} \to H{\bf Z}/p{\bf Z} \to H{\bf Z_p ^{\wedge}}$,
where ${\bf Z}$ is the integers and ${\bf Z_p ^{\wedge}}$ is the $p$-completion
of ${\bf Z}$, i.e. the ring of $p$-adic numbers. In other words, in the above setup,
we have taken $R$, $A$ and $B$ to be $H{\bf Z}$, $H{\bf Z}/p{\bf Z}$ and 
$H{\bf Z_p ^{\wedge}}$
respectively. First observe that $B \wedge_A B \simeq B$. Indeed, by Theorem IV.2.1 of 
~\cite{EKMM}, we have that
$$\pi_\ast( H{\bf Z_p ^{\wedge}} \wedge_{H{\bf Z}/p{\bf Z}} H{\bf Z_p ^{\wedge}}) \cong
Tor_\ast ^{{\bf Z}/p{\bf Z}} ({\bf Z_p ^{\wedge}}, {\bf Z_p ^{\wedge}}).$$
Hence, $\pi_0(B \wedge_A B) = 
\pi_0( H{\bf Z_p ^{\wedge}} \wedge_{H{\bf Z}/p{\bf Z}} H{\bf Z_p ^{\wedge}}) \cong
{\bf Z_p ^{\wedge}} \otimes_{{\bf Z}/p{\bf Z}} {\bf Z_p ^{\wedge}}$. 
However, by Theorem 7.2 of ~\cite{Eisen}, 
${\bf Z_p ^{\wedge}} \otimes_{\bf Z}/p{\bf Z} {\bf Z_p ^{\wedge}}$ is 
isomorphic to the $p$-adic completion of ${\bf Z_p ^{\wedge}}$, and since 
${\bf Z_p ^{\wedge}}$ is already complete, we conclude that 
$\pi_0(B \wedge_A B) \cong {\bf Z_p ^{\wedge}}$. As for $\pi_\ast(B \wedge_A B)$ for 
$\ast >0$, they are all trivial, since by Theorem 7.2 of ~\cite{Eisen}, 
${\bf Z_p ^{\wedge}}$ is flat over ${\bf Z}/p{\bf Z}$, and hence 
$Tor_\ast ^{\bf Z}/p{\bf Z} ({\bf Z_p ^{\wedge}}, {\bf Z_p ^{\wedge}}) \cong 0$ 
for $\ast >0$. Thus, we conclude that  $B \wedge_A B \simeq B$.

To prove that  $A \wedge_R B$ and $B \wedge_R B$ are not weakly equivalent, it is enough 
to show that $\pi_0(A \wedge_R B)$ is not isomorphic to $\pi_0(B \wedge_R B)$, which is 
evident, since 
$$\pi_0(A \wedge_R B) \cong ({\bf Z}/p{\bf Z}) \otimes_{\bf Z} {\bf Z_p ^{\wedge}} 
\cong {\bf Z_p ^{\wedge}}/p {\bf Z_p ^{\wedge}} \cong {\bf Z}/p{\bf Z},$$
while 
$\pi_0(B \wedge_R B) \cong {\bf Z_p ^{\wedge}} \otimes_{\bf Z} {\bf Z_p ^{\wedge}}$.

\end{example} 

To produce a sufficient condition for (~\ref{eq:general}) to hold, first we set up 
the notation, then introduce a few key identities which, if true, would imply 
the equation (~\ref{eq:general}). We discuss conditions under which these identities 
hold and, to conclude the section, summarize our findings in two (\'etale descent) 
lemmas.   

The objective is to compare the algebras $B \otimes_R X$ and 
$(A \otimes_R X) \wedge_A B$. We do this by comparing two towers of objects that 
approximate $B \otimes_R X$ and $(A \otimes_R X) \wedge_A B$ respectively. For the 
special case $X=S^1$ such towers were considered in ~\cite{Randy} for the category 
of chain complexes and adopted to the category of $S$-algebras in ~\cite{Minas}.  

Fix a simplicial set $X$. Define $I_A$ to be the hofiber of the multiplication map
$A \otimes_R X \to A$. Then $I_A$ inherits a multiplicative structure and we define 
$I_A/I_A ^n$ by the pushout diagram
$$
\xymatrix{
I_A ^n \ar[r] \ar[d]
&
I_A \ar[d]\\
\ast \ar[r]
&
I_A/I_A ^n,
}
$$
where the smash powers of $I_A$ are taken over $A \otimes_R X$. 

\begin{proposition}
\label{prop:towers}
Let $A \to B$ be thh-\'etale. Then the towers $\{(I_A/I_A ^n) \wedge_A B\}$ and 
$\{ I_B/I_B ^n \}$ are weakly equivalent. Consequently,
$$holim [(I_A/I_A ^n) \wedge_A B] \simeq  holim I_B/I_B ^n$$
\end{proposition}

\begin{proof}
We begin by showing that 
\begin{equation}
\label{eq:firststep}
I_B/I_B ^2 \simeq I_A/I_A ^2 \wedge_A B. 
\end{equation}

To see this, we employ new notation to denote the fiber of $A \otimes_R X \to A$ by 
$I_X$ whenever we wish to consider it as a functor of simplicial sets, as opposed to 
$R$-algebras. Observe that since $I_X/I_X ^2$ is a linear functor and 
$X \cong S^0 \wedge X$, we have that $I_X/I_X ^2$ is equivalent to 
$I_{S^0}/I_{S^0} ^2 \wedge X$. Recall that $I_{S^0}/I_{S^0} ^2 \simeq \Sigma TAQ(A|R)$
(see e.g. ~\cite{Minas}). Thus, to show (\ref{eq:firststep}), it suffices to prove 
that 
$$\Sigma TAQ(B|R) \wedge X \simeq \Sigma TAQ(A|R) \wedge_A B \wedge X,$$
which, in turn, is an immediate consequence of the transfer sequence of $TAQ$:
$$TAQ(A|R) \wedge_A B \to TAQ(B|R) \to TAQ(B|A),$$
combined with the fact that $TAQ(B|A) \simeq \ast$ since $A \to B$ is thh-\'etale. 
 
To complete the proof, we induct on $n$. Suppose the natural map 
$I_A/I_A ^{n-1} \wedge_A B \to I_B/I_B ^{n-1}$ is a weak equivalence. By naturality,
we have a commutative diagram

$$
\xymatrix{
(I_A ^{n-1}/I_A ^n) \wedge_A B \ar[r] \ar[d]
&
(I_A/I_A ^{n}) \wedge_A B \ar[r] \ar[d]
&
(I_A/I_A ^{n-1}) \wedge_A B \ar[d]\\
I_B ^{n-1}/I_B ^n \ar[r]
&
I_B/I_B ^n \ar[r]
&
I_B/I_B ^{n-1}
},
$$
where the objects in the left column are the hofibers of the right maps. Since both 
rows are (co)fibration sequences and the right vertical map is a weak equivalence by 
inductive assumption, it's enough to show that the left vertical map is also a weak 
equivalence. This, however, is an immediate consequence of Proposition 2.4 of 
~\cite{Minas}, which states that 
$$ [{\bigwedge ^n}_A I_A/I_A^2 ]_{h\Sigma_n} \simeq   I_A^n/I_A^{n+1}, $$
where the lower script $A$ in the above smash product on the left 
indicates that the smash product is taken over $A$. Thus, we have a series 
equivalences
$$ (I_A^n/I_A^{n+1})\wedge_A B \simeq  
[{\bigwedge ^n}_A I_A/I_A^2 ]_{h\Sigma_n} \wedge_A B \simeq
[{\bigwedge ^n}_B I_A/I_A^2 \wedge_A B]_{h\Sigma_n} \simeq
[{\bigwedge ^n}_B I_B/I_B^2 ]_{h\Sigma_n} \simeq
I_B^n/I_B^{n+1},
$$
which proves that the left vertical arrow, and consequently the middle one, are weak 
equivalences. 

\end{proof}

Observe that, in view of the above proposition, the \'etale descent formula 
(~\ref{eq:general}) will hold if $A \otimes_R X \wedge_A B$ and $B \otimes_R X$
are equivalent to $holim [(I_A/I_A ^n) \wedge_A B]$ and $holim [(I_B/I_B ^n)$ 
respectively. To address this, we pause to discuss completions and complete objects 
in our framework.  

\begin{definition}
{\bf 1.} Let $A$ be a cofibrant $R$-algebra. Define the completion $(A\otimes_R X)^\wedge $ 
of $A\otimes_R X$ to be the inverse limit $holim (A\otimes_R X)/I_A ^n$.

{\bf 2.} For an $A \otimes_R X$-module $M$, the completion $M^\wedge$ of $M$ is defined to be
$holim(M/I_A^n)$, where $M/I_A^n$ is the cofiber of the obvious map
$I_A ^n \wedge_{A\otimes_R X} M \to (A\otimes_R X) \wedge_{A\otimes_R X} M
\stackrel{\cong}{\to} M$. Here as before the powers of $I_A$ are taken over $A\otimes_R X$.

{\bf 3.} $M$ is complete if the natural map 
$M \to  holim[M/I_A ^n]$ is a weak equivalence.

\end{definition}

The following result helps to transmit information between an $S$-algebra and its 
completion.

\begin{proposition}
\label{prop:compl}
If $M$ is a finite $A$-CW-complex, then the natural map 
$$(A\otimes_R X)^\wedge \wedge_{A\otimes_R X} M\otimes_R X \to  (M\otimes_R X)^\wedge$$
is an equivalence. 

Consequently, if $B$ is a thh-\'etale algebra over $A$, which is a finite $A$-CW-complex
when viewed as an $A$-module,  
then the 
completion of $B \otimes_R X$ with respect to $I_B$ is weakly equivalent to the 
completion of $B \otimes_R X$ viewed as a $A \otimes_R X$-module. 
\end{proposition}

\begin{proof}
The Proposition is clearly true for $M=A$. Observe that as a consequence of adjunctions
$$\C_R(A \otimes_R X, B) \cong \U(X, \C_R(A,B)) \cong \C_R(A,F(X_+,B)),$$
we have that $(\Sigma^i A) \otimes_R X \cong \Sigma^i (A \otimes_R X)$, where $\Sigma^i A$
is the $i'th$ suspension of $A$; and in the above adjunctions $\C_R$ and $\U$ are the 
categories of commutative $R$-algebras and unbased spaces respectively. Hence, 
$((\Sigma^i A) \otimes_R X)^\wedge \simeq \Sigma^i (A \otimes_R X)^\wedge$, i.e. the 
Proposition holds for suspensions of $A$ as well. 

Now suppose, the statement is true for some module $K$ and let $F$ be a wedge of sphere 
modules $S_A ^i$ with a hocofiber $N$:
\begin{equation}
\label{eq:cw}
F \to K \to N.
\end{equation} 
Consider the following commutative diagram
$$
\xymatrix{
(A \otimes_R X)^\wedge \wedge_{A \otimes_R X} F \ar[r] \ar[d]
&
(A \otimes_R X)^\wedge \wedge_{A \otimes_R X} K \ar[r] \ar[d]
&
(A \otimes_R X)^\wedge \wedge_{A \otimes_R X} N \ar[d]\\
F^\wedge \ar[r]
&
K^\wedge \ar[r]
&
N^\wedge
}
$$
Note that both rows are cofibrations and the two left vertical maps are weak 
equivalences - the first one by our above discussion on suspensions of $A$, and 
the second one by assumption on K. Hence the right vertical map 
$(A \otimes_R X)^\wedge \wedge_{A \otimes_R X} N \to N^\wedge$ is also a 
weak equivalence, which proves the first part of the proposition, as $A$-CW-complexes
are built precisely via sequences (~\ref{eq:cw}). 

To prove the second part of the Proposition, we apply 
$- \wedge_{A \otimes_R X} (B \otimes_R X)$ to the sequence 
$I_A \to A \otimes_R X \to A$ to get a cofibration sequence
$$I_A\wedge_{A \otimes_R X} (B \otimes_R X) 
\to (A \otimes_R X) \wedge_{A \otimes_R X} (B \otimes_R X) 
\to A\wedge_{A \otimes_R X} (B \otimes_R X).$$
Note that by the base change formula (~\ref{eq:base}) for tensor products,
the last term $A \wedge_{A \otimes_R X} B \otimes_R X$ is equivalent to $B \tensor_A X$,
which, in turn, is weakly equivalent to $B$ by thh-\'etale assumption. Hence we have a 
cofibration sequence
$$I_A\wedge_{A \otimes_R X} (B \otimes_R X) \to B \otimes_R X \to B,$$
and are, thus, entitled to conclude that 
$ I_A\wedge_{A \otimes_R X} (B \otimes_R X) \simeq I_B$. The conclusion follows from the 
first part of the Proposition. 

\end{proof}

We are ready to state our first \'etale descent lemma. 

\begin{lemma}
\label{etdescomp}
(\'etale descent, complete case) Let $A$ be a cofibrant $R$-algebra, such that
$A\otimes_R X$ is complete, and $B$ be a cofibrant $A$-algebra which is a finite 
$A$-CW-complex when viewed as an $A$-module. Then $A \to B$ is 
thh-\'etale if and only if the \'etale descent formula holds:
$$  
(A \otimes_R X) \wedge_A B \simeq B \otimes_R X.
$$
\end{lemma}

\begin{proof}

We only need to prove the `only if' direction. Let $R \to A \to B$ be as in 
lemma, with $A \to B$ thh-\'etale. Then, by definition of completeness and 
due to the fact that smashing with finite CW-complexes commutes with holims, we 
have
\begin{equation}
\label{eq:compl1}
(A \otimes_R X) \wedge_A B \simeq 
holim[(A \otimes_R X)/I_A ^n] \wedge_A B \simeq
holim[((A \otimes_R X)/I_A ^n) \wedge_A B]
\end{equation}
Recall that by Proposition ~\ref{prop:towers}, 
\begin{equation}
\label{eq:compl2}
holim[((A \otimes_R X)/I_A ^n) \wedge_A B] \simeq
holim[(B \otimes_R X)/I_B ^n].
\end{equation}
Hence, it remains to prove that $holim[(B \otimes_R X)/I_B ^n]$ is weakly equivalent to 
$B \otimes_R X$, or in other words, that $B \otimes_R X$ is complete with respect to 
$I_B$, which, of course, is equivalent to being complete as an $A \otimes_R X$-module 
by Proposition ~\ref{prop:compl}. Denote the homotopy fiber of the natural map 
$A \otimes_R X \to ( A \otimes_R X)^\wedge$ by $K$ and consider the following diagram
whose right column is obtained by applying $- \wedge_{A \otimes_R X} (B \otimes_R X)$
to the cofiber sequence $K \to A \otimes_R X \to ( A \otimes_R X)^\wedge$:
$$
\xymatrix{
K \wedge_A B \ar[r] \ar[d]
&
K \wedge_{A \otimes_R X} (B \otimes_R X) \ar[d]\\
(A \otimes_R X) \wedge_A B \ar[r] \ar[d]
&
(A \otimes_R X) \wedge_{A \otimes_R X} (B \otimes_R X) \ar[d]
\cong B \otimes_R X \\
( A \otimes_R X)^\wedge \wedge_A B \ar[r]
&
( A \otimes_R X)^\wedge \wedge_{A \otimes_R X} (B \otimes_R X)
\simeq ( B \otimes_R X)^\wedge
}
$$
Since $A\otimes_R X$ is complete, $K$ is contractible; hence the top row is a weak 
equivalence. The bottom row is also an equivalence since combining equations 
(~\ref{eq:compl1}) and (~\ref{eq:compl2}) we get 
$$( A \otimes_R X)^\wedge \wedge_A B \simeq 
( A \otimes_R X) \wedge_A B \simeq
holim[(B \otimes_R X)/I_B ^n] \simeq ( B \otimes_R X)^\wedge .$$
Hence, we are allowed to conclude that the middle row is also an equivalence, which 
proves the lemma.  

\end{proof}

\begin{remark}
We would like to point out that it is this \'etale 
descent lemma that prompted us to consider the thh-\'etale algebras (in addition to
\'etale ones). Of course, the more direct translation of the `\'etale' notion from 
discrete algebra appears to be what we have defined as \'etale $S$-algebras,
since in both cases \'etale essentially means unramified, 
i.e with a vanishing module of differentials. Hence, perhaps one would like/hope
to prove an \'etale descent lemma with an \'etale condition (as opposed to a slightly 
stronger thh-\'etale requirement as we have imposed). However, as we have demonstrated,
the (stronger) thh-\'etale condition is a necessary one. We also note that the notion 
of thh-\'etale maps is also a generalization of \'etale maps from discrete algebra;
in fact, as pointed out earlier, when restricted to Eilenberg-MacLane spectra \'etale
and thh-\'etale coincide.   

\end{remark}

We return to the completeness assumption in the \'etale descent lemma above. 
That assumption is satisfied if $A$ is connective and the simplicial set $X$ is 
such that  $\pi_0(X)=0$, as clearly the connectivity of maps 
$$A \otimes_R X \to (A \otimes_R X)/I_A ^n$$
increases with $n$, since with $A$ connective and $X$ connected, $I_A$ is at least 
1-connected. Equivalently, the connectivity of fibers $I_A ^n/I_A ^{n+1}$ 
increases with $n$. Moreover, if $B$ is a connective $A$-algebra then by 
Eilenberg-Moore spectral sequence (Section 4, Chapter IV of ~\cite{EKMM}),
the connectivity of the maps
$$A \otimes_R X \wedge_A B \to ((A \otimes_R X)/I_A ^n) \wedge_A B $$   
also increases with $n$, which implies that 
$$A \otimes_R X \wedge_A B \simeq holim[((A \otimes_R X)/I_A ^n) \wedge_A B]. $$
By Proposition~\ref{prop:compl}, $ holim[((A \otimes_R X)/I_A ^n) \wedge_A B]$ is 
weakly equivalent to $ holim[(B \otimes_R X)/I_B ^n]$, which, in turn is equivalent
to $B \otimes_R X$ since $B$ is connective and $X$ is connected, and hence, 
$B \otimes_R X$ is complete. 

We have proved the following lemma.

\begin{lemma} 
\label{etdesconn}
(\'etale descent, connective case)
Let $A$ be a connective cofibrant $R$-algebra, $B$ a connective cofibrant 
$A$-algebra, and $X$ a connected simplicial set. Then $A \to B$ is thh-\'etale
if and only if the \'etale descent formula folds:
$$(A \otimes_R X) \wedge_A B \simeq B \otimes_R X.$$
\end{lemma}

In conclusion of this section, we present a result that helps to detect the 
condition $B \otimes_A X \simeq B$ necessary (and often sufficient) for the 
\'etale descent  
Lemmas~\ref{etdescomp} and ~\ref{etdesconn} to hold. We set up the notation first. 

For a simplicial set $X_\ast$, let $J_X$ be the fiber of the obvious (induced by 
multiplication) map $B \otimes_A X \to B$ to emphasize that $J$ is a functor of 
simplicial sets. 

\begin{proposition}
Let $A \to B$ be a map of commutative $R$-algebras and $X_\ast$ a simplicial set 
such that $B \otimes_A X$ is complete with respect to $J_X$. 
Then $B \otimes_A X \simeq B$ if and only if 
$H_\ast(X, TAQ(B|A))=0$ for all $\ast$. 
\end{proposition}

\begin{proof}

We begin by observing that $B \otimes_A X \simeq B$ if and only if $J_X \simeq \ast$.
This in turn implies that $ J_X /J_X ^2 \simeq \ast$. Furthermore, the converse 
of this is also true. Indeed, let $ J_X /J_X ^2 \simeq \ast$. By ~\cite{Randy} or 
~\cite{Minas} we have that 
$$hofiber(J_X /J_X ^{n+1} \to J_X /J_X ^n) \simeq 
[(J_X /J_X ^2)^{\wedge n}]_{h\Sigma_n}.$$
This result is listed as Proposition 2.4 in ~\cite{Minas}, which in turn is the 
adaptation to the framework of $S$-algebras of a similar result obtained in 
~\cite{Randy} for the category of chain complexes.  Now if 
$ J_X /J_X ^2 \simeq \ast$ then the first term and
all homotopy fibers in the inverse limit system
$\{J_X /J_X ^n\}$ are contractible. Hence, 
$J_X \simeq holim J_X /J_X ^n \simeq \ast$. 

Recalling that the term $ J_X /J_X ^2$ is linear and that $X \cong X \wedge S^0$,
we get an identity $J_X /J_X ^2 \simeq X \wedge J_{S^0}/J_{S^0} ^2$. Thus,
$B \otimes_A X \simeq B$ if and only if $X \wedge J_{S^0}/J_{S^0} ^2 \cong \ast$.
To complete the proof it remains to observe that 
$TAQ(B|A) \simeq J_{S^0}/J_{S^0} ^2$, and hence $B \otimes_A X \simeq B$ is 
equivalent to $X \wedge TAQ(B|A) \simeq \ast $, or in other words, to 
$H_\ast(X, TAQ(B|A))=0$ for all $\ast$. 

\end{proof}  

\section{HKR Theorem}

\begin{theorem}
\label{theorem1}
Let $f:R \to A$ be thh-smooth in the category of connective 
$S$-algebras. Then the natural (derivative) map 
$THH(A|R) \to \Sigma TAQ (A|R)$ has a section in the category of $A$-modules 
which induces an equivalence of $A$-algebras:
$$\P_A \Sigma TAQ (A|R) \stackrel{\simeq}{\longrightarrow} THH(A|R).$$
\end{theorem} 

\begin{proof}
First we show that the Theorem holds for polynomial extensions 
$R \rightarrow \P_R X$, where $X$ is a cell $R$-module. 
Our first objective is to compute $TAQ(\P_R X|R)$. While one can do this directly 
from definitions, we present a somewhat more concise computation that employs series 
of adjunctions. By Proposition 3.2 of ~\cite{Maria}, for every $\P_R X$-module $M$,
$$ h \M_{\P_R X} (TAQ(\P_R X|R), M) \cong h\C_{R/\P_R X}(\P_R X,\P_R X \vee M),$$
where $\C_{R/\P_R X}$ is the category of $R$-algebras over $\P_R X$, and 
$h \M$ and $h \C$ indicate the corresponding homotopy categories. Of course, it 
is immediate that 
$\C_{R/\P_R X}(\P_R X,\P_R X \vee M) \cong \C_R(\P_R X, M)$. Furthermore, since the 
free functors $\P_R$ and $ \P_R X \wedge_R -$ (with $X$ a cell $R$-module) are left 
adjoints which preserve cofibrations and trivial cofibrations, 
they induce adjunctions on homotopy categories as well (see ~\cite{Dwyer}). 
Thus, we get
$$ h\C_R(\P_R X, M) \cong h \M_R(X, M) \cong h M_{\P_R X}(\P_R X \wedge_R X, M).$$
Hence, by Yoneda's lemma,  we have an equivalence of $\P_R X$ modules  
$TAQ(\P_R X|R) \simeq \P_R X \wedge_R X$. 

On the other hand, by a theorem of McClure, Schw\"anzl and Vogt(~\cite{Vogt}), 
$THH(\P_R X|R) \cong \P_R X \otimes_R S_\ast ^1$. We have adjunction homeomorphisms
$$\C_R(\P_R X \otimes_R S_\ast ^1, B) \cong \U(S_\ast ^1, \C_R(\P_R X, B))
\cong \U(S^1, \M_R(X, B)) \cong \M_R(X \wedge S^1 _+, B) \cong 
\C_R(\P_R (X \wedge S^1 _+), B),$$
where $\C_R$ is the category of commutative $R$-algebras, $\U$ is the category 
of unbased topological spaces, and $B$ is a commutative $R$-algebra. Hence, by 
Yoneda's lemma, 
$THH(\P_R X|R) \cong \P_R (X \wedge S^1 _+)$ as $R$-algebras. 
Of course, $\P_R (X \wedge S^1 _+)$ (and consequently $THH(\P_R X|R)$) also has a 
structure of a $\P_R X$-algebra, which is more evident once we observe that
$\P_R (X \wedge S^1 _+) \cong \P_R (X \vee \Sigma X) \cong \P_R X \wedge_R 
\P_R (\Sigma X)$. Finally, note that by the base change formula for polynomial 
algebras, we have 
$$\P_R (X \wedge S^1 _+) \cong \P_R X \wedge_R \P_R (\Sigma X) 
\cong \P_{\P_R X}(\P_R X \wedge_R \Sigma X).$$
 



\noindent
Hence, recalling that $TAQ(\P_R X|R) \simeq \P_R X \wedge_R X$,
we are allowed to conclude that as $\P_R X$-algebras 
Topological Hochschild Homology $THH(\P_R X|R)$ is equivalent to
$\P_{\P_R X} (\P_R X \wedge_R \Sigma X) \cong  \P_{\P_R X} (\Sigma TAQ(\P_R X|R))$.
 
Now let $R \to A$ be an arbitrary smooth map. Thus we have a family of sequences 
$$ R \to \P_R X \stackrel{\phi}{\rightarrow} A_\alpha$$
with $\phi$ thh-\'etale. By ~\cite{Maria}, this sequences give rise to cofibration sequences
\begin{equation}
TAQ(\P_R X|R) \wedge_{\P_R X} A_\alpha \rightarrow
TAQ(A_\alpha|R) \rightarrow TAQ(A_\alpha|\P_R X)
\end{equation}
Since $\phi$ is thh-\'etale, the last term of this sequence is 0. Hence,
\begin{equation}
\label{eq:taqcov}
TAQ(\P_R X|R) \wedge_{\P_R X} A_\alpha \stackrel{\simeq}{\rightarrow}
TAQ(A_\alpha|R).
\end{equation}
Similarly, the sequences $R \to A \to A_\alpha$ produce cofibration sequences
\begin{equation}
TAQ(A|R) \wedge_{A} A_\alpha \rightarrow
TAQ(A_\alpha|R) \rightarrow TAQ(A_\alpha|A)
\end{equation}
Since the maps $A \to A_\alpha$ are thh-\'etale by definition, $TAQ(A_\alpha|A)$ 
are contractible. Hence, we get an equivalence of $A$-modules
\begin{equation}
\label{eq:taqar}
TAQ(A|R) \wedge_{A} A_\alpha \stackrel{\simeq}{\rightarrow}
TAQ(A_\alpha|R).
\end{equation}
Combining the above Lemma~\ref{etdesconn} with the fact that we have proved the 
theorem for polynomial extensions, we get a series of equivalences
\begin{equation}
THH(A_\alpha|R) \cong
THH(\P_R X|R) \wedge_{\P_R X} A_\alpha \cong
\P_{\P_R X} (\Sigma TAQ(\P_R X|R)) \wedge_{\P_R X} A_\alpha
\end{equation}
Next, observe that 
$\P_{\P_R X} (\Sigma TAQ(\P_R X|R)) \wedge_{\P_R X} A_\alpha \cong
\P_{A_\alpha} (\Sigma TAQ(\P_R X|R) \wedge_{\P_R X} A_\alpha)$, which 
combined with the equation (~\ref{eq:taqcov}) gives us the theorem for the 
extensions $R \to A_\alpha$:
\begin{equation}
\label{eq:pretheor}
THH(A_\alpha|R) \simeq  \P_{A_\alpha} (\Sigma TAQ(A_\alpha|R)).
\end{equation}
To complete the proof, note that the Lemma~\ref{etdesconn} applied to the 
thh-\'etale map $A \to A_\alpha$ gives an equivalence 
$THH(A_\alpha|R) \simeq THH(A|R) \wedge _A A_\alpha$; 
and plugging this and the 
equation (~\ref{eq:taqar}) into the above equivalence (~\ref{eq:pretheor}), we get
$$THH(A|R) \wedge _A A_\alpha \simeq 
\P_{A_\alpha} (\Sigma TAQ(A|R) \wedge_{A} A_\alpha) \simeq
\P_{A} (\Sigma TAQ(A|R)) \wedge_{A} A_\alpha.
$$
 
Recalling the second condition of the definition of thh-\'etale covers 
$A \to A_\alpha$, we conclude that $THH(A|R)$ and 
$\P_{A} (\Sigma TAQ(A|R))$ are equivalent as $A$-algebras. 

\end{proof}

\begin{theorem}
Let $B \to R \stackrel{f}{\to} A$ be maps of connective $S$-algebras with 
$f$ thh-smooth. Then the first fundamental sequence of modules of 
differentials splits, i.e
$$TAQ(A|B) \simeq (TAQ(R|B) \wedge_R A) \vee TAQ(A|R)$$
\end{theorem}

\begin{proof}

Consider the suspension of the first fundamental sequence of differential modules
for $B \to R \stackrel{f}{\to} A$:
\begin{equation}
\label{eq:fund}
\Sigma TAQ(R|B) \wedge_R A \to \Sigma TAQ(A|B) \to \Sigma TAQ(A|R).
\end{equation}
By Theorem ~\ref{theorem1}, we have a map $\Sigma TAQ(A|R) \to THH(A|R)$ which 
is a section to the derivative map. The smash product over $B$ of the maps
$id:A \to A$ and $B \to R$ induces a map 
$THH(A|R) \to THH(A \wedge_B R|R) \stackrel{\simeq}{\to} THH(A|B)$. 
Thus we get a map 
$$\phi: \Sigma TAQ(A|R) \to THH(A|B).$$
Next consider the natural commutative diagram
$$
\xymatrix{
THH(A|B) \ar[r] \ar[d]
&
THH(A|R) \ar[d]\\
\Sigma TAQ(A|B) \ar[r]
&
\Sigma TAQ(A|R)
}
$$
It's easy to see that the map $\phi: \Sigma TAQ(A|R) \to THH(A|B)$ is 
a section to the map from $THH(A|B)$ to  $\Sigma TAQ(A|R)$ in the above diagram. 
Thus, $\phi$ composed with the derivative map $THH(A|B) \to \Sigma TAQ(A|B)$ gives a map 
$\Sigma TAQ(A|R) \to \Sigma TAQ(A|B)$ which is a section to the second map in 
the first fundamental sequence, and thus splits the sequence. This map combined with 
the first map in the fundamental sequence (~\ref{eq:fund}) induces a map 
\begin{equation}
\label{eq:split}
(TAQ(R|B) \wedge_R A) \vee TAQ(A|R) \to TAQ(A|S).
\end{equation}
Since the second map in the equation (~\ref{eq:fund}) has a section, it is surjective 
on homotopy groups and the long exact sequence of homotopy groups associated to the 
cofibration sequence (~\ref{eq:fund}) breaks up into a series of {\it{split}} short exact 
sequences:
$$ \pi_i(TAQ(R|B) \wedge_R A) \to \pi_i(TAQ(A|S)) \to \pi_i (TAQ(A|R)).$$
Hence, $\pi_i(TAQ(A|S)) \cong \pi_i(TAQ(R|B) \wedge_R A) \oplus \pi_i (TAQ(A|R))$, 
which implies that the map (~\ref{eq:split}) induces an isomorphism on homotopy groups and
is thus a weak equivalence. 

\end{proof}

\end{document}